\documentclass[11pt]{amsart}

\usepackage{amssymb} 
\usepackage[mathscr]{eucal} 
\usepackage[all]{xy}
\usepackage{amsmath} 

\def\cal#1{{ \mathcal {#1} }} 

\newcommand{\TS}{\textstyle}

\newcommand{\tagnum}[1] 
{\noindent \hbox to 22pt {\hss\bf #1.}}

\theoremstyle{plain}
\newtheorem{bigthm}{Theorem}

\newtheorem{thm}{Theorem}[section]

\newtheorem{lem}[thm]{Lemma}
\newtheorem{prop}[thm]{Proposition}

\newtheorem{rem}[thm]{Remark}

\numberwithin{equation}{section}

\def\rxar#1{\xrightarrow{\kern -1pt #1}}
\def\lxar#1{\xleftarrow{#1 \kern -3pt}}

\newcommand{\Pn}[1]{{ \mathbb{P}^{#1}}\relax}

\newcommand{\ra}{\ensuremath{\rightarrow}}

\def\eea{\end{eqnarray*}}
\def\bea{\begin{eqnarray*}}

\def\F{{\mathbb{F}}}

\def\PP{{\mathbb{P}}}

\newcommand{\shExt}{\cal Ext\kern 2pt}
\newcommand{\shHom}{\cal Hom\kern 2pt}
\newcommand{\Syz}{\cal Syz\kern 2pt}

\DeclareMathOperator{\Supp}{Supp}

\title{On Wahl's proof of $\mu(6)=65$}

\author{Roberto Pignatelli}

\address{Roberto Pignatelli, Dipartimento di Matematica Univ. di Trento\\
Via Sommarive, 14 38050-Trento (Italy)}

\author{Fabio Tonoli}

\address{Fabio Tonoli, Dipartimento di Matematica Univ. di Trento\\
Via Sommarive, 14 38050-Trento (Italy)}

\email{roberto.pignatelli@unitn.it}
\email{tonoli@science.unitn.it}


\begin{document}

\maketitle

\section*{Introduction}

In this note we present a short proof of the following theorem of
D. Jaffe and D. Ruberman: 

\smallskip\noindent\textbf{Theorem \cite{JR}.}
{\em A sextic hypersurface in $\PP^3$ has at most 65 nodes.}
\smallskip

The bound is sharp by Barth's construction \cite{Ba} of a sextic
with $65$ nodes.

Following Beauville \cite{Bea}, to a set of $n$ nodes on a surface is
associated a linear subspace of $\F^n$ (where $\F$ is the field with
two elements) whose elements corresponds to the so-called 
{\em even subsets} of the set of the nodes. Studying this code
Beauville proved that the maximal number of nodes of a quintic surface
is $31$.

The same idea was used by Jaffe and Ruberman, but
their proof is not so short as the one of Beauville, partly because at
that time a complete understanding of the possible cardinalities of an
even set of nodes was missing.

Almost at the same time, J. Wahl \cite{wahl} proposed a much shorter
proof of the same result. He proved indeed the following (see the
beginning of the next section for the missing definitions)

\smallskip\noindent\textbf{Theorem \cite{wahl}.}
{\em Let $V\subset \F^{66}$ be a code, with weights in
$\{24,32,40\}$. Then $\dim(V)\leq12$.}
\smallskip

He claimed that Jaffe-Ruberman's theorem follows as a corollary since
the code associated to a nodal sextic has dimension at least $n-53$  
(see section $1$ of \cite{CaTo} for this computation).
In fact, he used an incorrect result stated by Casnati and Catanese in
\cite{CaCa}, asserting that the possible
cardinalities of an even set of nodes on a sextic are only $24,32$ and
$40$. 
Recently Catanese and Tonoli showed indeed

\smallskip\noindent\textbf{Theorem \cite{CaTo}.}
{\em  On a sextic nodal surface in $\Pn 3$, an even set 
of nodes has cardinality in $\{24,32,40,56\}$.}
\smallskip

\noindent
Note however that \cite{CaTo} used a result by Jaffe and Ruberman, 
namely that there is no even set of nodes of cardinality $48$.

By the above theorem the proof of the theorem of Jaffe and Ruberman
reduces to the following
\begin{bigthm}
Let $V\subset \F^{66}$ be a code with weights in
$\{24,32,40,56\}$. Then $\dim(V)\leq12$.
\end{bigthm}

This statement is in fact theorem 8.1 of \cite{JR}. 
Anyway, its proof is much more complicated than Wahl's one and
moreover requires computers computations. 
In this short note we give an elementary proof, 
using and integrating Wahl's ideas.

{\bf Acknowledgement:} Both authors would like to thank F. Catanese for
suggesting this problem and for the reading of an earlier version of
this paper. His comments have been very helpful in clarifying some arguments.

\section{Notation and general results from coding theory}

A {\bf code} is (in this note) a vector subspace $V\subset\F^n$, where
$\F$ is the field with two elements. 
A {\bf word} is a vector $v=(v_1,\ldots,v_n)\in \F^n$. Its {\bf support} $\Supp(v)$
is the set $\{i \mid v_i\neq 0\}$ 
of coordinates that do not vanish in $v$, 
its {\bf weight} $|v|$ is the cardinality of its support. 
The {\bf length} of a code is
the cardinality of the union of the supports of all its elements. 
A code $V\subset\F^n$ is said to be {\bf spanning} if it has length $n$.

A code is {\bf even} if all its words have even weight, {\bf doubly
  even} if all its weights are divisible by $4$. The number of words
of weight $i$ in the code $V$ is denoted by $a_i(V)$ or simply
$a_i$ when no confusion arises. The \textbf{weight enumerator} of the
code $V$ is the homogeneous polynomial
$$W_V(x,y)=\sum a_i x^{n-i}y^i.$$

The standard scalar product in $\F^n$ associates to each code  
its dual code , {\it i.e.}, its annihilator $V^*\subset \F^n$, 
which has complementary dimension.
We set $a_i^*:=a_i(V^*)$.

\begin{rem}\label{cazzatine}
1) $V\subset\F^n$ is spanning if and only if $a_1^*=0$.

2) If $v^* \in V^*$ has weight 2,  
the subset of $V$ given by all words $v$ with $\Supp(v)\cap \Supp(v^*) =
\emptyset$ is a subcode of codimension at most $1$ (and length at
most $n-2$).

3) A doubly even code is automatically isotropic, {\it i.e.}, $V
\subset V^*$.
\end{rem}

The \textbf{MacWilliams identity} (cf. \cite{mw}) 
states that the weight enumerator $W_{V^*}(x,y)$ of
the dual code $V^*$ equals $W_V(x+y,x-y)/2^d$, {\it i.e.},
\begin{equation}  \label{MW}
  \sum a_i^* x^{n-i} y^i = 
\frac{1}{2^d} \left(\sum a_i (x+y)^{n-i} (x-y)^i\right).
\end{equation}

As explained in \cite{wahl}, comparing the
coefficients of $x^{n-i}y^{i}$ for $i\leq 3$ in both sides 
of (\ref{MW}) gives (since $a_0=a_0^*=1$):
\begin{lem}\cite[Lemma 2.4]{wahl}\label{lem2.4}
Let $V\subset \F^n$ be a spanning code of dimension $d$. Then:
\begin{subequations}
\begin{gather}
\TS \sum_{i>0} a_i = 2^d-1\label{1}\\
\TS \sum i a_i = 2^{d-1}n\label{2}\\
\TS \sum i^2 a_i = 2^{d-1}(a_2^*+n(n+1)/2)\label{3}\\
\TS \sum i^3 a_i = 2^{d-2}\left(3(a_2^*n-a_3^*)+n^2(n+3)/2\right)\label{4}
\end{gather}
\end{subequations}
\end{lem}

The following proposition gives dimension and weights 
of a projected linear code.

\begin{prop}\cite[Prop. 2.8]{wahl}\label{prop2.8}
Let $V\subset\F^n$ be a code of dimension $d$. 
Fix a word $w\in V$ and consider the projection 
$\pi\colon  \F^n \ra \F^{n-|w|}$ onto the complement of the support
of $w$. Then
\begin{enumerate}
\item If $w$ is not a sum of two disjoint words in $V$, then 
$V':=\pi(V)$ is a code of dimension $d'=d-1$.
\item $|\pi(v)|=\frac12(|v|+|v+w|-|w|)$. 
\end{enumerate}
\end{prop}

\begin{proof}
If $\ker \pi_{|V}$ contains, besides $w$, another
word $v$, one can write a disjoint sum $w=v+(w-v)$. 
Thus, in the hypothesis of (1), $\dim \ker
\pi_{|V}=1$ and therefore $d'=d-1$.

For $(2)$, let $r$ be the cardinality of the intersection of the two
supports of $v$ and $w$. Then $|v|=r+|\pi(v)|$ and $|v|+|w|=|v+w|+2r$.
\end{proof}

\section{The proof}

\begin{lem}\cite[Lemma 2.6]{wahl}\label{lem2.6}
The dimension of a code with weights in $\{24, 32\}$ is at most $9$.
\end{lem}
\begin{proof} Let $n$ be the length of the code and $d$ its dimension.
Solving the linear system given by (\ref{1}) and (\ref{2}),
$a_{24}=2^{d-4}(64-n)-4$, $a_{32}=2^{d-4}(n-48)+3$. Substituting in
(\ref{3})
$$
2^8\left(2^{d-6}\cdot 9 \cdot (2^6-n) + 2^{d-2} \cdot (n-48)+3\right)=2^{d-1}(a_2^*+n(n+1)/2)
$$ 
If $d>9$, then $2^{d-1}$ divides the R.H.S. but not the L.H.S., a contradiction.
\end{proof}

\begin{rem}\label{rem1}
A code $V\subset\F^{67}$ with weights $\geq 24$ has necessarily
$a_{56}\leq 1$.
\end{rem}

\begin{proof}
Indeed, if there are two different words of weight 56, their sum
has weight at least 24 and then the cardinality of the intersection of
their supports is at least $1/2(56+56-24)=44$. Therefore their span
has length  $\geq 44+2\cdot(56-44)=68$.
\end{proof}

\begin{lem}\label{lem2}
The dimension of a code  $V\subset \F^{67}$ with weights in $\{24,
32, 56\}$ is at most $10$. 
\end{lem}

\begin{proof}
If $a_{56}=0$ the result follows by Lemma \ref{lem2.6}. 

Otherwise, by Remark \ref{rem1}, $a_{56}=1$.
The intersection of $V$ with any hyperplane not containing its unique
word of weight 56 is a code $V'$ of dimension $\dim(V)-1$ with weights in
$\{24, 32 \}$ and the result follows again by Lemma \ref{lem2.6}.
\end{proof}


%
%
%
%

\begin{proof}[Proof of Theorem A]
Suppose that there exists a code $V\subset
\F^{66}$ with  weights in $\{24,32,40,56\}$ of dimension $13$. 
Let $n$ be its length and consider $V$ as a spanning code in $\F^n$.

By Lemma \ref{lem2} we have $a_{40}>0$.
For each word $w\in V$ with weight 40 we consider the projection
$\pi_w$ onto the complement of the support of $w$.
By Proposition \ref{prop2.8}, $V':=\pi_w(V)\subset \F^{n-40}$ is a
doubly even code of dimension $12$.
So $V'$ is an isotropic subspace, $n-40\geq 24$ and we obtain
$n\geq 64$: more precisely $n \in \{64,65,66\}$. 
\smallskip

Suppose $n=64$. 
For each word $w\in V$ of weight $40$, $\pi_w(V)$ is isotropic of
dimension $12$ in $\F^{24}$, so $\pi_w(V)=(\pi_w(V))^*$. Let $\mathbb
I\in \F^{24}$ be the vector with all coordinates 1: $\mathbb I\in (\pi_w(V))^*$
(since $\pi_w(V)$ is even) and therefore $\mathbb I\in\pi_w(V)$.

If $v\in V$ is a word such that both the weights $|v|,|v+w|$ are 
$\leq 40$, then by Proposition \ref{prop2.8} $|\pi_w(v)|\leq 20$; therefore
by remark \ref{rem1} $a_{56}(V)=1$ and $\mathbb I=\pi_w(\overline{v})$ for the
unique word $\overline{v} \in V$ with $|\overline{v}|=56$. 

Fix one coordinate not in the support of $\overline{v}$ and let
$V'' \subset V$ be the subcode defined by the vanishing of the given
coordinate. Since $\mathbb I=\pi_w(\overline{v})$, the support of $w$
contains the complementary of the support of $\overline{v}$: then $w
\not\in V''$. Since this holds for each $w \in V$ with $|w|=40$, then
$V''$ has no word of weight $40$: it is a code
of dimension $12$ with weights in $\{24,32,56\}$, contradicting lemma
\ref{lem2}.

\smallskip

Suppose $n=65$.
Solving the equations (\ref{1})-(\ref{4}), we obtain 
$a_{56}=\frac12(a_2^*-a_3^*-5)$ and thus $a_2^*> 0$. Let then $z \in
V^*$ be a word of length $2$.

For each word $w\in V$ of weight $40$, $a_2^*(\pi_w(V))=0$: in fact,
for any word $z'\in(\pi_w(V))^*$ of weight 2, ${\rm Span}(V',z')$
is an isotropic subspace of dimension 13 in $\F^{25}$, absurd.
Therefore every word $w$ of weight 40 satisfies $\Supp(w)\supset \Supp(z)$.

By remark \ref{cazzatine} the subset of $V$ given by all words $v$
with $\Supp(v)\cap \Supp(z) = \emptyset$ is a subcode of dimension at
least $12$ with weights in $\{24,32,56\}$, contradicting Lemma \ref{lem2}.
\smallskip

Then $n=66$.
Solving the equations (\ref{1})-(\ref{4}), we obtain 
$a_{56}=a_2^*-\frac12(a_3^*+13)$ and thus $a_2^*\geq7$. We choose two
words $z_1\neq z_2$ in $V^*$ of weight 2.

If we show that for each word $w\in V$ of weight $40$,
$a_2^*(\pi_w(V))\leq 1$, then $\Supp(w)$ intersects $Z=\Supp(z_1)\cup
\Supp(z_2)$. Therefore, by remark \ref{cazzatine}, the subset of $V$
given by all words $v$ with $\Supp(v)\cap Z = \emptyset$ is
a code of dimension at least 11 and weights among$\{24,32,56\}$,
contradicting again Lemma \ref{lem2}. 

\smallskip

So it remains to show only that for each word $w\in V$ of weight $40$,
$a_2^*(\pi_w(V))\leq 1$.

If $z'\in (\pi_w(V))^*$ is a word of weight 2, then
$V'':={\rm Span}(\pi_w(V),z')\subset \F^{26}$ is an isotropic subspace of
dimension 13, and thus $\mathbb I\in V''=(V'')^*$. 
Being $\pi_w(V)$ doubly even, $\mathbb I,z'\in V''\setminus \pi_w(V)$,
and therefore $\mathbb I+z'$ is a word in $\pi_w(V)$ of weight 24. 
Thus $a_2^*(\pi_w(V)) \leq a_{24}(\pi_w(V))$.

If $v\in V$ is a word such that both the weights $|v|,|v+w|$ are 
$\leq 40$, then by Proposition \ref{prop2.8} $|\pi_w(v)|\leq 20$; therefore
$a_{24}(\pi_w(V))\leq a_{56}(V)\leq 1$ (the last inequality by remark
\ref{rem1}). 

\end{proof}





\begin{thebibliography}{[McW-Sl}

\bibitem[Ba]{Ba}
W.~Barth,
\newblock {Two projective surfaces with many nodes, admitting the
symmetries of the icosahedron.}
\newblock {{\em J. Alg. Geom.}, 5 (1996), 173--186 .}

\bibitem[Be]{Bea}
A.~Beauville,
\newblock {Sur le nombre maximum de point doubles d'une surface dans
$\Pn 3$ ($\mu(5)=31$)}.
\newblock {In {\em Algebraic Geometry, Angers, 1979,
A.~Beauville ed., Sijthoff \& Noordhoff}, (1980), 207--215.}
%
%
%
%
%
%
%
%
\bibitem[Ca-Ca]{CaCa}
G. Casnati and F. Catanese,
\newblock {Even sets of nodes are bundle symmetric.}
\newblock {{\em J. Diff. Geom.}, 47 (1997), 237--256.}
\newblock {
Erratum.}
\newblock {{\em J. Diff. Geom.}  50  (1998),  no. 2, 415.}



\bibitem[Ca-To]{CaTo}
F. Catanese, F. Tonoli,
\newblock {Even sets of nodes on sextic surfaces}
\newblock {To appear in {\em J.E.M.S}, 2006.}

%
%
%
%
%
%
%
%
%
\bibitem[Ja-Ru]{JR}
D.~B.~Jaffe and D.~Ruberman,
\newblock {A sextic surface cannot have 66 nodes.}
\newblock {{\em J. Alg. Geom.}, 6 (1997), 151--168 .}

%
%
%
\bibitem[McW-Sl]{mw}
F. J. MacWilliams, N. J. A. Sloane,
\newblock {The theory of error-correcting codes I, II.}
\newblock{{\em North-Holland Mathematical Library}, Vol. 16.
North-Holland Publishing  Co., Amsterdam-New York-Oxford, (1977),
           i--xv + 1--369, resp.  i--ix + 370--762.}

%
%
%
%
%
%
%
%
\bibitem[Wa]{wahl}
J. Wahl,
\newblock {Nodes on sextic hypersurfaces in $P\sp 3$}.
\newblock {\em J. Differential Geom.} 48 (1998), no. 3, 439--444.

%
%
%
%
%



\end{thebibliography}
\end{document}